\documentclass[12pt]{article}
\usepackage{amssymb}
\usepackage{amsmath}

\begin{document}

\title{On exit times of L\'evy-driven Ornstein--Uhlenbeck processes}

\author{Konstantin Borovkov\footnote{Department of Mathematics and Statistics, The
University of Melbourne, Parkville 3010, Australia. Supported by ARC Centre of
Excellence for Mathematics and Statistics of Complex Systems.} \ and Alexander
Novikov\footnote{Department of Mathematical Sciences, University of Technology, PO Box
123, Broadway, Sydney, NSW 2007, Australia. Supported by an ARC Discovery grant.}}

    \date{}

    \maketitle

\begin{abstract}We prove two martingale identities which involve exit
times of L\'evy-driven Ornstein--Uhlenbeck processes. Using these identities we find an
explicit formula for the Laplace transform of the exit time under the assumption that
positive jumps of the L\'evy process are exponentially distributed.
\end{abstract}



\section{Introduction}

Let $X_{t}$, $t\geq 0$, be an Ornstein--Uhlenbeck (O-U) process driven by a L\'evy
process $L_{t}$, i.e. $X_{t}$ is a solution of the equation
\begin{equation*}
X_{t}=x-\beta \int_{0}^{t}X_{s}ds+L_{t},~~\,t\geq 0.
\end{equation*}

We assume that the parameter\thinspace $\beta $ is positive and the initial
value $X_{0}=x$ is non-random.

In the special case when $L_{t}$ is a compound Poisson process, the process $%
X_{t}$ is also known in applications as a \textquotedblleft
shot-noise\textquotedblright\ process or a \textquotedblleft storage
process\textquotedblright\ with a linear release function.

One of the most important for the models of that sort problems is to
determine or to approximate the distribution of the first passage time
\begin{equation*}
\tau _{b}=\inf \{t>0:\,X_{t}\geq b\}
\end{equation*}%
of a given level $b>x$. The problem was discussed for the Gaussian
O-U processes by Darling \& Siegert (1953). Explicit
representations for the Laplace transform $E(e^{-\mu \tau _{b}})$
were found in the papers of Hadjiev (1983) and Novikov (1990,
2003) in the case when $L_{t}$ has no positive jumps (the
so-called spectrally negative case). Moreover, the papers of
Novikov and \`{E}rgashev (1993) and Novikov (2003) provide some
bounds and asymptotic approximations for the distribution of $\tau
_{b}$. In particular, it was proved in Novikov (2003), Theorem 2,
that the distribution of $\tau _{b}$ is exponentially bounded
under the condition that $L_{t}$ has a diffusion part or positive
jumps. The papers Perry et al. (2001) and Borovkov and Novikov
(2001) contain some general results on integral equations for the
distributions of $\tau_b$ and $X_{\tau_b}.$

It seems that the first results for L\'evy-driven O-U processes with exponentially
distributed jumps were obtained by Tsurui and Osaki (1976) for the case when the
parameter $1/\beta $ is an integer. For the case of arbitrary $\beta >0$ and
exponentially distributed positive jumps, explicit formulas for the Laplace transform
$E(e^{-\mu \tau _{b}})$ and the expectation of $\tau _{b}$ can be found in the paper of
Novikov et al. (2005) who solved the corresponding integro-differential equation.
Recently, Jacobsen and Jensen (2006) found the joint Laplace transform $E(e^{-\mu \tau
_{b}+wX_{\tau _{b}}})$ in the form of a linear combination of contour integrals under
the assumption that a distribution of the positive jumps is a mixture of exponential
ones.

In what follows,  we always assume that the following condition holds:
\begin{equation}
E\log (1+|L_{1}|)<\infty   \label{wolfe}
\end{equation}%
(this is a sufficient and necessary condition for convergence of $X_{t}$ in
distribution to a proper limit, see e.g. Wolfe (1982)).

In Section 2 we prove two martingale identities (see Theorems 1 and  2
below) which involve both the first passage time $\tau _{b}$ and $X_{\tau
_{b}}.$ These identities enable one to obtain explicit bounds for the
distribution of $\tau _{b}$ (e.g. an explicit lower bound for $E\tau _{b}$)
just by neglecting the overshoot
\begin{equation*}
\chi _{b}=X_{\tau _{b}}-b.
\end{equation*}

In Section 3 we use Theorems 1 and 2 for deriving explicit representations
of the Laplace transform and the mean of $\tau _{b}$ under the assumption
that the positive jumps of $L_{t}$ are exponentially distributed but without
any restrictions on the distribution of the negative jumps of $L_{t}.$ We
also prove the Exponential Limit Theorem for $\tau _{b}$ as $b\rightarrow
\infty $.

\section{Martingale identities}

In what follows we always assume that $L_t$ has a non-zero component with positive
jumps(or, equivalently, $\Pi (0,\infty)>0$ where  $\Pi (dx)$ is the L\'evy-Khinchin
measure associated  with $L_{t}$). This assumption implies that $L_{t}$ has the
following representation:
\begin{equation}
L_{t}=Q_{t}+R_{t}  \label{represL}
\end{equation}%
with the compound Poisson process
\begin{equation}
R_{t}=\sum_{k=1}^{N_{t}}\xi _{k},  \label{represR}
\end{equation}%
where $\xi _{k}$ are the jumps of $L_{t}$ which are greater than
some positive number $\delta ,\newline P\{\xi _{k}>\delta \}>0,$
$N_{t}$ is a Poisson process with rate $E(N_{1})=\lambda >0.$ We
also assume that the componet $Q_{t}$ can only contain  a
diffusion part and jumps less than or equal to $\delta $ and
therefore $Q_t$ and $R_t$ are independent.

Set
\begin{equation*}
K=\sup \{u\geq 0:Ee^{uL_{1}}<\infty \}.
\end{equation*}%
We shall further assume that $K>0$ and set
\begin{equation}
\varphi (u)=\frac{1}{\beta }\int_{0}^{u}\frac{\log (Ee^{vL_{1}})}{v}%
dv,~\,0\leq u<K.  \label{phi}
\end{equation}%
Since $Q_{t}$ and $R_{t}$ are independent, we have
\begin{equation}
\varphi (u)=\Delta (u)+W(u),  \label{deltaplus}
\end{equation}%
where we put%
\begin{equation}
\Delta (u)=\frac{1}{\beta }\int_{0}^{u}\frac{\log (Ee^{vQ_{1}})}{v}dv,
\label{delta}
\end{equation}%
\begin{equation*}
\text{ }W(u)=\frac{1}{\beta }\int_{0}^{u}\frac{\log (Ee^{vR_{1}})}{v}dv=%
\frac{\lambda }{\beta }\int_{0}^{u}\frac{(Ee^{v\xi _{1}}-1)}{v}dv.
\end{equation*}%
Under the assumption (\ref{wolfe}) the integrals in (\ref{phi}) and (\ref%
{delta}) converge (see some details of the proof for this fact in Wolfe
(1982) or Novikov (2003)) and so $\varphi (u),~\Delta (u)$ and $W(u)$ are
finite continuous functions. Besides, for all $u\in \lbrack 0,\infty )$ the
following lower bound holds\footnote{$c$ and $C$ are some positive constants}
\begin{equation}
\Delta (u)\geq -c-Cu  \label{lowDelta}
\end{equation}
(see Novikov (2003)).

Using the inequality $e^{x}>1+x+x^{2}/2,\,x>0,$ we obtain also
\begin{equation}
\text{ }W(u)\geq \frac{\lambda }{\beta }(u\delta +u^{2}\delta ^{2}/4)P\{\xi
_{k}>\delta \}>0.  \label{lowW}
\end{equation}%
Set%
\begin{equation*}
\varphi (K)=\lim_{u\uparrow K}\varphi (u).
\end{equation*}%
If $K=\infty ,$ then $\varphi (K)=\infty $ due to (\ref{lowW}). If $%
0<K<\infty ,$~then the value $\varphi (K)$ could be finite or infinite as
illustrated by the following example where the Compound Poisson process $%
L_{t}=\sum_{k=1}^{N_{t}}\xi _{k},~$has the jumps $\xi _{k}$ with the Gamma
distribution, i.e.%
\begin{equation*}
P(\xi _{k}\in dx)=\frac{x^{\rho -1}e^{-x}}{\Gamma (\rho )}dx,~x>0,\rho >0.
\end{equation*}%
Then $K=1$ and by direct calculation
\begin{equation*}
\varphi (u)=\frac{\lambda }{\beta }\int_{0}^{u}\frac{1-(1-v)^{\rho }}{%
(1-v)^{\rho }v}dv,~~u<1,
\end{equation*}%
so that%
\begin{equation*}
\varphi (1)<\infty \text{ ~~}\emph{for}\mathit{~}~\rho <1
\end{equation*}%
and
\begin{equation*}
\varphi (1)=\infty \text{ ~~\emph{for} ~}\rho \geq 1.
\end{equation*}%
Note that when $\rho =1$ (this is the case of exponentially distributed
jumps with mean one) we have the explicit formula
\begin{equation}
\varphi (u)=-\frac{\lambda }{\beta }\log (1-u),~~u<1.  \label{expjump}
\end{equation}

Set
\begin{equation}
G(z,\mu )=\int\limits_{0}^{K}e^{uz-\varphi (u)}u^{\mu -1}du,~~~\mu >0.
\label{G(z)}
\end{equation}%
This function is, obviously, finite when $K<\infty .~$\ For the case $%
K=\infty $ the finiteness of $G(z,\mu )$ is implied by (\ref{lowDelta}) and (%
\ref{lowW}).

\medskip\noindent
\textbf{Theorem 1. }\emph{Let condition }(\ref{wolfe})\emph{\ hold,~}$%
0<K\leq \infty $ \emph{and }$\varphi (K)=\infty .\,$\emph{Then}
\begin{equation}
E(e^{-\mu \beta \tau _{b}}G(X_{\tau _{b}},\mu ))=G(x,\mu ),\,\,\,\,\,\mu >0.
\label{mar3}
\end{equation}

\medskip\noindent
\textbf{Proof.} First consider the case%
\begin{equation*}
K=\infty ,
\end{equation*}%
in which it was shown by Novikov (2003) that the process $e^{-\mu \beta
t}G(X_{t},\mu )$ is a martingale with respect to the natural filtration $%
\mathcal{F}_{t}=\sigma \{L_{s},s\leq t\}.$

Applying the optional stopping theorem, we have for any $t>0$%
\begin{equation}
E[I\{\tau _{b}\leq t\}e^{-\mu \beta \tau _{b}}G(X_{\tau _{b}},\mu
)]+E[I\{\tau _{b}>t\}e^{-\mu \beta t}G(X_{t},\mu )]=G(x,\mu ).  \label{mar31}
\end{equation}%
Since $X_{t}\leq b$ on the event $\{\tau _{b}>t\}$ and $G(x,\mu )$ is a
nondecreasing function of $x,$ we~obtain
\begin{equation*}
E[I\{\tau _{b}>t\}e^{-\mu \beta t}G(X_{t},\mu )]\leq e^{-\mu \beta t}P\{\tau
_{b}>t\}G(b,\mu )\rightarrow 0
\end{equation*}%
as $t\rightarrow \infty $. The first term on the LHS of (\ref{mar31})
clearly converges monotonically to $E(e^{-\mu \beta \tau _{b}}G(X_{\tau
_{b}},\mu ))$ as $t\rightarrow \infty $ because $\tau _{b}$ is finite with
probability one (in fact, it is even exponentially bounded). So (\ref{mar3})
holds when $K=\infty .$

To prove (\ref{mar3}) in the case $0<K<\infty ,$ we shall truncate positive
jumps of $L_{t}$ by a positive constant $A$ and then justify a passage to
the limit as $A\rightarrow \infty .$

Set%
\begin{equation*}
L_{t}^{A}=Q_{t}+R_{t}^{A}
\end{equation*}%
with
\begin{equation*}
R_{t}^{A}=\sum_{k=1}^{N_{t}}\min (\xi _{k},A).
\end{equation*}%
Let $X_{t}^{A}$ be an O-U process driven by $L_{t}^{A}$,
\begin{equation*}
\tau _{b}^{A}=\inf \{t\geq 0:\,X_{t}^{A}\geq b\}
\end{equation*}%
and
\begin{equation*}
\varphi _{A}(u)=\frac{1}{\beta }\int_{0}^{u}\frac{\log (Ee^{vL_{1}^{A}})}{v}%
dv=\Delta (u)+W(u,A),
\end{equation*}%
where we put%
\begin{equation}
W(u,A)=\frac{\lambda }{\beta }\int_{0}^{u}\frac{(Ee^{v\min (\xi _{1},A)}-1)}{%
v}dv.  \label{preseWA}
\end{equation}%
It is obvious from the L\'evy-Khintchin formula that the right distribution tail of
$L_{1}^{A}$ decays faster than any exponential function, so that the respective value
$K=K(A)=\infty .$ Hence, identity (\ref{mar3}) does hold for the process $X_{t}^{A}:$
\begin{equation}
E\left[ e^{-\mu \beta \tau _{b}^{A}}\int\limits_{0}^{\infty }\exp \{uX_{\tau
_{b}^{A}}^{A}-\varphi _{A}(u)\}u^{\mu -1}du\right] =\int\limits_{0}^{\infty
}e^{ux-\varphi _{A}(u)}u^{\mu -1}du.  \label{mar5}
\end{equation}%
Further we note~that as $\,A\rightarrow \infty $
\begin{equation}
\int\limits_{0}^{\infty }e^{ux-\varphi _{A}(u)}u^{\mu -1}du\rightarrow
\int\limits_{0}^{K}e^{ux-\varphi (u)}u^{\mu -1}du  \label{mar6}
\end{equation}%
which gives the RHS of (\ref{mar3}). To see this, we note that, as $%
A\rightarrow \infty ,$ for any $u<K$%
\begin{equation*}
W(u,A)\rightarrow W(u),\ ~~\varphi _{A}(u)\rightarrow \varphi (u),\,
\end{equation*}%
and for $u\geq K$
\begin{equation*}
W(u,A)\rightarrow \infty .
\end{equation*}%
and, oviously, the last two relations imply (\ref{mar6}). \bigskip Next, on
the LHS of (\ref{mar5}) we write%
\begin{equation*}
\int\limits_{0}^{\infty }=\int\limits_{0}^{K}+\int\limits_{K}^{\infty }
\end{equation*}%
and consider convergence of the corresponding two terms separately. Note
that
\begin{equation}
X_{\tau _{b}}^{A}\leq b+\delta +\min (\xi _{N_{\tau _{b}^{A}}},A).
\label{oversh}
\end{equation}%
Obviously, $\tau _{b}^{A}$ could only decrease as $A$ increases.
Choose now a positive constant $A_{0}$ such that $P\{\xi
_{1}<A_{0}\}>0$ and so $\tau _{b}^{A_{0}}$ is exponentially
bounded. Then we have for all $A>A_0$
\begin{equation}
e^{uX_{\tau _{b}}^{A}}\leq e^{u(b+\delta )}\sum_{k=1}^{N_{\tau
_{b}^{A_{0}}}}e^{u\min (\xi _{k},A)} \label{oversh2}
\end{equation}
where by Wald's identity
\begin{equation}
E(\sum_{k=1}^{N_{\tau _{b}^{A_{0}}}}e^{u\min (\xi
_{k},A)})=E(N_{\tau _{b}^{A_{0}}})E(e^{u\min (\xi _{1},A)})
\label{Wald}
\end{equation}%
and
\begin{equation*}
E(N_{\tau _{b}^{A_{0}}})=\lambda E(\tau _{b}^{A_{0}})<\infty .
\end{equation*}%
Collecting together the above bounds we obtain for the $\int\limits_{K}^{%
\infty }-$part of the LHS of (\ref{mar5}) the following bound:
\begin{multline*}
E\left[ e^{-\mu \beta \tau _{b}^{A}}\int\limits_{K}^{\infty }\exp \{uX_{\tau
_{b}^{A}}^{A}-\varphi _{A}(u)\}u^{\mu -1}du\right] \\
\leq C\int\limits_{K}^{\infty
}E(e^{u\min (\xi _{1},A)})e^{u(b+\delta )-\varphi _{A}(u)}u^{\mu -1}du,
\end{multline*}%
where $C=\lambda E(\tau _{b}^{A_{0}}).~$To show that the last integral
converges to zero as $A\rightarrow \infty ,$ we note that (\ref{preseWA})
implies
\begin{equation*}
E(e^{u\min (\xi _{1},A)})=\frac{\beta u}{\lambda }\frac{\partial W(u,A)}{%
\partial u}+1.
\end{equation*}%
This means that
\begin{multline}
\int\limits_{K}^{\infty }E(e^{u\min (\xi _{1},A)})e^{u(b+\delta )-\varphi
_{A}(u)}u^{\mu -1}du \\
=\frac{\beta }{\lambda }\int\limits_{K}^{\infty }\frac{\partial W(u,A)}{%
\partial u}e^{u(b+\delta )-W(u,A)-\Delta (u)}u^{\mu
}du+\int\limits_{K}^{\infty }e^{u(b+\delta )-W(u,A)-\Delta (u)}u^{\mu -1}du.
\label{mar7}
\end{multline}%
The last integral tends to zero as $A\rightarrow \infty $ due to the fact that
$W(u,A)\rightarrow \infty $ for $u\geq K.$

Integrating by parts the first integral on the RHS of (\ref{mar7}), we
obtain:
\begin{equation*}
\int\limits_{K}^{\infty }\frac{\partial W(u,A)}{\partial u}e^{u(b+\delta
)-W(u,A)-\Delta (u)}u^{\mu }du=-\int\limits_{K}^{\infty }e^{u(b+\delta
)-\Delta (u)}u^{\mu }d(e^{-W(u,A)})
\end{equation*}%
\begin{equation*}
=e^{K(b+\delta )-\Delta (K)}K^{\mu }e^{-W(K,A)}+\int\limits_{K}^{\infty
}e^{-W(u,A)}d(e^{u(b+\delta )-\Delta (u)}u^{\mu }).
\end{equation*}%
Now it should be clear that the last two terms converge  to zero as $%
A\rightarrow \infty $ due to the fact that $W(u,A)\rightarrow \infty $ for $%
u\geq K.$ \ So, we have proved the convergence of the $\int\limits_{K}^{%
\infty }-$part of the LHS of (\ref{mar5}) to zero.

To study the part with $\int\limits_{0}^{K}$ on the LHS of (\ref{mar5}),
note that the random variable $X_{\tau _{b}^{A}}^{A}$ coincides with $%
X_{\tau _{b}}$ on the set $\{\max\limits_{k\leq N_{\tau _{b}^{A}}}\xi
_{k}<A\}$ (because no jumps are truncated up to the time $\tau _{b}\leq \tau
_{b}^{A}$). Obviously, as $A\rightarrow \infty $%
\begin{equation}
P\{\max\limits_{k\leq N_{\tau _{b}^{A}}}\xi _{k}<A\}\rightarrow 1,
\label{prob1}
\end{equation}%
and hence
\begin{multline*}
E\left[ I\{\max\limits_{k\leq N_{\tau _{b}^{A}}}\xi _{k}<A\}e^{-\mu \beta \tau
_{b}^{A}}\int\limits_{0}^{K}e^{uX_{\tau _{b}}-\varphi _{A}(u)}u^{\mu -1}du\right]
 \\
\rightarrow E\left[ e^{-\mu \beta \tau _{b}}\int\limits_{0}^{K}e^{uX_{\tau
_{b}}-\varphi (u)}u^{\mu -1}du\right] .
\end{multline*}%
To complete the proof, we need to check only that%
\begin{equation}
\lim_{A\rightarrow \infty }E\left[ I\{\max\limits_{k\leq N_{\tau
_{b}^{A}}}(\xi _{k})\geq A\}\int\limits_{0}^{K}e^{uX_{\tau
_{b}^{A}}^{A}-\varphi _{A}(u)}u^{\mu -1}du\right] =0.  \label{kk}
\end{equation}%
To see this, note that in view of (\ref{oversh2}), we have for all $A\geq
A_{0}$
\begin{equation*}
\int\limits_{0}^{K}e^{uX_{\tau _{b}^{A}}^{A}-\varphi _{A}(u)}u^{\mu
-1}du\leq \sum_{k=1}^{N_{\tau _{b}^{A_{0}}}}\int\limits_{0}^{K}e^{u\min (\xi
_{k},A)}e^{u(b+\delta )-\Delta (u)}e^{-W(u,A)}u^{\mu -1}du\leq C\eta _{A},
\end{equation*}%
where we put%
\begin{equation*}
\eta _{A}=\sum_{k=1}^{N_{\tau _{b}^{A_{0}}}}\int\limits_{0}^{K}e^{u\min (\xi
_{k},A)}e^{-W(u,A)}u^{\mu -1}du,~C=e^{K(b+\delta )-\min_{u\leq K}\Delta (u)}.
\end{equation*}%
Due to this bound and (\ref{prob1}), for the validity of (\ref{kk}) it is
sufficient to show that $\{\eta _{A},A\geq A_{0}\}$ is a family of uniformly
integrable \ random variables or, equivalently, that%
\begin{equation*}
\lim_{A\rightarrow \infty }E(\eta _{A})=E(\lim_{A\rightarrow \infty }\eta
_{A})<\infty .
\end{equation*}%
In view of (\ref{Wald}) the latter property is equivalent to%
\begin{equation*}
\lim_{A\rightarrow \infty }E(\int\limits_{0}^{K}e^{u\min (\xi
_{k},A)}e^{-W(u,A)}u^{\mu -1}du)=E(\int\limits_{0}^{K}e^{u\xi
_{k}}e^{-W(u)}u^{\mu -1}du)<\infty .
\end{equation*}%
which one can easily verify (e.g. using monotonicity of the functions $\min (\xi
_{k},A)$ and $W(u,A)).$

This completes the proof.

\medskip\noindent
\textbf{Theorem 2. }\emph{Let condition (\ref{wolfe})\ hold,~}$0<K<\infty $%
\emph{\ and}$~\varphi (K)=\infty .$ \emph{Then}%
\begin{equation*}
\beta E(\tau _{b})=E\int\limits_{0}^{K}(e^{uX_{\tau
_{b}}}-e^{ux})e^{-\varphi (u)}u^{-1}du.
\end{equation*}

\medskip\noindent
\textbf{Proof. }We will derive this identity from Theorem 1 by passing to
the limit as $\mu \rightarrow 0.~\,$\ To justify this procedure we observe
that (\ref{mar3}) can be written in the following form:
\begin{equation*}
E(e^{-\mu \beta \tau _{b}}\int\limits_{0}^{K}(e^{uX_{\tau
_{b}}}-e^{ux})e^{-\varphi (u)}u^{\mu -1}du)=(1-Ee^{-\mu \beta \tau
_{b}})\int\limits_{0}^{K}e^{ux-\varphi (u)}u^{\mu -1}du,\,\,\,\mu >0.
\end{equation*}%
Here the LHS converges to $E\int\limits_{0}^{K}(e^{uX_{\tau
_{b}}}-e^{ux})e^{-\varphi (u)}u^{-1}du$ as $\mu \rightarrow 0$. One can
easily see (e.g. using integration by parts) that
\begin{equation}
\int\limits_{0}^{K}e^{ux-\varphi (u)}u^{\mu -1}du\thicksim \frac{1}{\mu }.
\label{exp4a}
\end{equation}%
This implies that
\begin{equation*}
\lim_{\mu \rightarrow 0}(1-Ee^{-\mu \beta \tau
_{b}})\int\limits_{0}^{K}e^{ux-\varphi (u)}u^{\mu -1}du=\beta E(\tau _{b})
\end{equation*}%
which completes the proof of Theorem 2.

\medskip\noindent
\textbf{Remark. }The assertion of Theorem 2 was proved for the case $%
K=\infty $ in Novikov (2003) under an additional assumption that%
\begin{equation*}
E(L_{1}^{-})^{\delta }<\infty \,\,\emph{for\,\,some\,\,}\delta >0.
\end{equation*}

\section{Exponentially distributed positive jumps}

In this section we use the same notation as in Section 2 and
assume that the process $Q_{t}$ in the decomposition
(\ref{represL}) does not contain positive jumps while the process
$R_{t}$ is a compound Poisson process with \ exponentially
distributed positive jumps, $E(\xi _{k})=K,~0<K<\infty ;~N_{t}$ is
a Poisson process with rate $E(N_{1})=\lambda >0.$ Note that under
these assumptions
\[e^{-\varphi (u)}=(1-u/K)^{\lambda /\beta}e^{-\Delta(u)}.
\]

\medskip\noindent
\textbf{Theorem 3. }\emph{For any }$\mu >0,$%
\begin{equation}
E(e^{-\mu \beta \tau
_{b}})=\frac{\int\limits_{0}^{K}(1-u/K)^{\lambda /\beta}u^{\mu
-1}e^{ux-\Delta (u)}du}{\int\limits_{0}^{K}(1-u/K)^{\lambda /\beta
-1}u^{\mu -1}e^{ub-\Delta (u)}du},  \label{exp1}
\end{equation}%
\begin{equation}
E(\tau _{b})=\frac{1}{\beta }\int%
\limits_{0}^{K}(e^{ub}-e^{ux}(1-u/K))(1-u/K)^{\lambda /\beta -1}e^{-\Delta
(u)}u^{-1}du.  \label{exp2}
\end{equation}%
\emph{Besides, as }$b\rightarrow \infty ,$%
\begin{equation}
E(\tau _{b})=Ce^{Kb}(Kb)^{-\lambda /\beta }(1+o(1)),~C=\frac{\Gamma (\lambda
/\beta )}{\beta K}e^{-\Delta (K)}  \label{exp5}
\end{equation}%
\emph{and the Exponential Limit Theorem holds:}
\begin{equation}
P\{\frac{\tau _{b}}{E(\tau _{b})}>x\}\rightarrow e^{-x}\,,~x>0.  \label{exp6}
\end{equation}

\medskip\noindent
\textbf{Proof. }Formulas (\ref{exp1}) and (\ref{exp2}) are direct
consequences of Theorem 1, Theorem 2 and the following two well-known facts
(which hold due to the memory-less property of the exponential distribution,
see a similar statement in Borovkov (1976) for the case $\beta =0)$:

1) the overshoot $\chi _{b}=X_{\tau _{b}}-b$ has the density
\begin{equation*}
p_{\chi _{b}}(x)=\frac{1}{K}e^{-x/K},~~~x>0;
\end{equation*}

2)
\begin{equation*}
\chi _{b}\,~~\,\emph{and}\text{ \ \ }\tau _{b}\text{ \ }\emph{%
are\,\,independent.}
\end{equation*}%

Relation (\ref{exp2}) implies, using the change of variables $(1-u/K)b=w,$
\begin{equation*}
\frac{d}{db}E(\tau _{b})=\frac{1}{\beta }\int\limits_{0}^{K}e^{ub}(1-u/K)^{%
\lambda /\beta -1}e^{-\Delta (u)}du
\end{equation*}%
\begin{equation*}
=\frac{1}{\beta }e^{Kb}b^{-\lambda /\beta
}\int\limits_{0}^{b}e^{-Kw}w^{\lambda /\beta -1}e^{-\Delta (K(1-w/b))}dw.
\end{equation*}%
Since the function $\Delta (K(1-w/b))$ is continuous and bounded in $w\in
\lbrack 0,K],$ the last integral converges as $b\rightarrow \infty $ to
\begin{equation*}
e^{-\Delta (K)}\int\limits_{0}^{\infty }e^{-Kw}w^{\lambda /\beta
-1}dw=e^{-\Delta (K)}\Gamma (\lambda /\beta )K^{-\lambda /\beta }.
\end{equation*}%
Hence
\begin{equation*}
\frac{d}{db}E(\tau _{b})\thicksim \frac{\Gamma (\lambda /\beta )}{\beta }%
e^{-\Delta (K)}e^{Kb}(bK)^{-\lambda /\beta },~b\rightarrow \infty .
\end{equation*}%
By well-know facts of theory of asymptotic expansions (see e.g. Olver
(1997)) it implies%
\begin{equation*}
E(\tau _{b})\thicksim \frac{\Gamma (\lambda /\beta )}{\beta K}e^{-\Delta
(K)}e^{Kb}(bK)^{-\lambda /\beta }
\end{equation*}%
and therefore we have proved (\ref{exp5}).

To derive (\ref{exp6}), we write the denominator in (\ref{exp1})
as follows:
\begin{equation}
\int\limits_{0}^{K}(e^{ub}-e^{ux}(1-\frac{u}{K}))(1-\frac{u}{K})^{\lambda
/\beta -1}u^{\mu -1}e^{-\Delta (u)}du+\int\limits_{0}^{K}(1-\frac{u}{K}%
)^{\lambda /\beta }u^{\mu -1}e^{ux-\Delta (u)}du . \label{exp6a}
\end{equation}%
Set $$\mu =z/(\beta E\tau _{b})$$ with a fixed $z>0$ and $E(\tau
_{b})$ defined in (\ref{exp2}). Clearly,$\mu \to 0$ as $
b\rightarrow \infty .$ Due to (\ref{exp4a}), the second term in
(\ref{exp6a}) (which is also the nominator in (\ref{exp1})) can
now be written as
\begin{equation}
\int\limits_{0}^{K}(1-u/K)^{\lambda /\beta }u^{\mu -1}e^{ux-\Delta (u)}du=%
\frac{\beta E(\tau _{b})}{z}(1+o(1)).\label{exp6ab}
\end{equation}%
Using (\ref{exp2}), the first term in (\ref{exp6a}) can be written
as
\begin{equation}
\int\limits_{0}^{K}(e^{ub}-e^{ux}(1-u/K))(1-u/K)^{\lambda /\beta -1}u^{\mu
-1}e^{-\Delta (u)}du=\beta E(\tau _{b})+\delta (b),  \label{exp6b}
\end{equation}%
where
\begin{equation*}
\delta (b)=\int\limits_{0}^{K}(e^{ub}-e^{ux}(1-u/K))(1-u/K)^{\lambda /\beta
-1}u^{-1}(u^{\mu }-1)e^{-\Delta (u)}du.
\end{equation*}%
Note that for $u>0$
\begin{equation*}
|u^{\mu }-1|=|e^{\mu \log u}-1|\leq \mu \max (u^{\mu },1)|\log u|.
\end{equation*}%
This implies~
\begin{equation*}
|\delta (b)|\leq \mu \max (K^{\mu },1)
\int\limits_{0}^{K}(e^{ub}-e^{ux}(1-u/K))(1-u/K)^{\lambda
/\beta-1}e^{-\Delta (u)}u^{-1}|\log u|du.
\end{equation*}%
Applying the same change of variables as above, one can easily
show that the last integral is  $O(e^{Kb}b^{-\lambda /\beta })$ as
$b\rightarrow \infty .$

Taking into account (\ref{exp5}), due to the setting for $\mu$ we
get \[\delta(b)=O(1).\] Now making the substitution in
(\ref{exp1}) relations (\ref{exp6a}), (\ref{exp6ab}),
(\ref{exp6b}) with the last result, we obtain
\begin{equation*}
E(e^{-z\tau _{b}/E(\tau _{b})})=\frac{\frac{\beta E(\tau _{b})}{z}(1+o(1))}{%
\beta E(\tau _{b})+O(1)+\frac{\beta E(\tau _{b})}{z}(1+o(1))}\rightarrow
\frac{1}{z+1}.
\end{equation*}%
Since the function $\frac{1}{1+z}$ is the Laplace transform of the
exponential distribution with mean $1,$ this completes the proof.

 \medskip
\noindent \textbf{Acknowledgement.} The authors are thankful to G.Miteteli and a
referee for constructive comments.


\begin{thebibliography}{99}

\bibitem{1} Borovkov, A. (1973) \emph{Stochastic processes in queueing theory},  Springer-Verlag, New York.


\bibitem{2} Borovkov, K., and Novikov, A. (2001)  {On a piece-wise deterministic Markov
process model.} {\em Statist. Probab. Lett.} 53, no. 4, 421--428.

\bibitem{3} Darling, D. A., and Siegert, A. J. F. (1953)
{The first passage problem for a continuous Markov process.} {\em Ann. Math.
Statistics} 24, 624--639.

\bibitem{4} Hadjiev, D. (1983)  {The first passage problem for
generalized Ornstein-Uhlenbeck processes with nonpositive jumps.} In: {\em S\'{e}%
minaire de probabilit\'{e}s,} XIX, 1983/84,  80--90, Lecture Notes in Math., 1123,
Springer, Berlin.

\bibitem{5} Jacobsen, M. and Jensen, A. (2006)  {\em Exit times for a class
of piecewise exponential Markov processes with two-sided jumps. } Dept. of
Applied Mathematics and Statistics, University of Copenhagen. Preprint No 5.

\bibitem{6} Kella, O. and Stadje, W. (2001)  {On hitting times for
compound Poisson dams with exponential jumps and linear release rate.} {\em J. Appl.
Prob.} 38, no. 3, 781--786.

\bibitem{7} Novikov, A. A. (1990)  {On the first exit time of an
autoregressive process beyond a level and an application to the `change-point' problem.
} {\em Theory Probab. Appl.} 35, no. 2, 269--279.

\bibitem{8} Novikov, A. A. and \`{E}rgashev, B. A. (1993).  {Limit
theorems for the time of crossing a level by an autoregressive process.} In: {\em Trudy
Mat. Inst. Steklova.} 202, 209--233. [In Russian. English translation in: {\em Proc.
Steklov Math. Inst.\/} 1994, 4 (202), 169--186.]

\bibitem{9} Novikov, A.A. (2003)  {Martingales and first-exit times for
the Ornstein--Uhlenbeck process with jumps.} {\em Theory Probab. Appl.}  48, 340--358.

\bibitem{10} Novikov, A.A., Melchers, R.E., Shinjikashvili, E. and
Kordzakhia, N. (2005)  {First passage time of filtered Poisson process with exponential
shape function}. {\em Probabilistic Engineering Mechanics,}   20, no. 1, 33-44.

\bibitem{11} Olver, F.W.J. (1997) \emph{Asymptotics and special functions.} 2nd end. AK Peters, Wellesley,
Mass.

\bibitem{12} Perry, D., Stadje, W. and Zacks, S. (2001) { First-exit times for Poisson shot
noise.} {\em Stoch. Models,}  17, no. 1, 25--37.

\bibitem{13} Tsurui, A. and Osaki, Sh. (1976)  {On a first-passage
problem for a cumulative process with exponential decay.} {\em Stochastic Processes
Appl.} 4, no. 1, 79--88.

\bibitem{14} Wolfe, S. (1982)  {On a continuous analogue of the
stochastic differential equation} $X_{n}=\rho X_{n-1}+B_{n}.\,$ {\em Stoch.\ Proc.\
Appl.,} 12, 301-312
\end{thebibliography}
\end{document}